\documentclass{article}
\usepackage{amsmath, amssymb, graphicx}
\title{Where Do the Terms of the Power Series Expansions of Sine and Cosine Functions Come from? Involutes!}
\author{{\bf V. N. Krishnachandran}\\ Vidya Academy of Science \& Technology\\ Thrissur -- 680501, Kerala, India \\ (email: {\tt krishnachandranvn@gmail.com})}
\date{}
\allowdisplaybreaks
\usepackage{tikz}
\usetikzlibrary{intersections, calc}
\newcommand\DrawControlAbove[3]{
node[#2,circle,fill=#2,inner sep=1pt,label={[black] above:{\normalsize#3}}] at #1 {}}
\newcommand\DrawControlLeft[3]{
node[#2,circle,fill=#2,inner sep=1pt,label={[black] left: {\normalsize#3}}] at #1 {}}
\newcommand\DrawControlRight[3]{
node[#2,circle,fill=#2,inner sep=1pt,label={[black] right:{\normalsize#3}}] at #1 {}}
\newcommand\ma{1.22173}
\newcommand\mb{0.34907}
\newcommand\mc{0.40724}
\newcommand\mr{9}
\newcommand\ms{4}
\newcommand\side{2.351141}
\begin{document}
\maketitle
\begin{abstract}
In the 1930's, a Russian school teacher Y. S. Chaikovsky presented a proof of the power series expansion of the sine and cosine functions without using calculus. In doing so he also showed the geometrical meanings of the various terms in these power series expansions. Chaikovsky's ideas were first published by Leo S. Gurin as a Note in the American Mathematical Monthly in 1996. The proofs, though they use only elementary mathematics, require some combinatorial arguments which may be hard even for bright students.

This paper is an attempt to bring the ideas of Chaikovsky once
 again to the attention of the mathematics teachers and students.
 The emphasis in this paper is on establishing the geometrical
 meanings of the various terms in the power series expansions of
 sine and cosine functions. This has been done taking recourse 
to methods of calculus. 

Since the concept of an involute is crucial in this investigation, we have  also included  a brief discussion on the concept of an involute of a curve and also a {\em Maxima} procedure for computing the parametric equations of a given curve. 
\end{abstract}
\section{Introduction}
It is well known that the definitions, terminology, concepts and
properties of the trigonometric sine and cosine functions evolved 
over a long period of time. In this long historical process Indian
mathematicians have made some profound contributions. The roots
of the terms `sine' and `cosine' have been traced to the Sanskrit
 words {\em jya} and {\em ko-jya}.  The construction of the 
first sine table has been attributed the the celebrated fifth-century CE Indian mathematician-astronomer Aryabhata \cite{Arya}. Much later,  
the following power series expansions of these functions were discovered by another Indian mathematician-astronomer  
Sangamagrama Madhava in the fifteenth century CE:
\begin{align}
\sin x & = x-\frac{x^3}{3!}+\frac{x^5}{5!}- \cdots \\
\cos x & = 1 - \frac{x^2}{2!} + \frac{x^4}{4!} - \cdots 
\end{align}
Madhava also had proofs for these expansions which are `satisfactory' even by modern standards of rigor. 
It was only two centuries later, after the discovery of calculus, that these power series expansions came to be known in Europe. In the Western world, it was Newton who first established these power series expansions. Now, the proofs of these power series expansions have reduced to simple exercises in applications of the so-called Taylor's theorem in calculus. 

But there is a certain mystery behind the series. For, in the functions $\sin x$ and $\cos x$ the quantity $x$ is an angle and  the power series expansions contain terms like $x^2, x^3, \ldots$. These terms represent multiplication of an angle with another angle. The process of multiplying two angles has apparently no geometrical meaning and hence the process is completely counter-intuitive. Of course an angle is a dimension-less quantity and so may be treated as pure number. However there is no geometrical process which produces the square of an angle. 

In the 1930's, a Russian school teacher Y. S. Chaikovsky presented before his bewildered students in the eighth grade a proof of the power series expansion of the sine and cosine functions without using calculus \cite{Gurin}.  His arguments brought to light the geometrical meanings of the various terms in the power series expansions of $\sin x$ and $\cos x$. Chaikovsky did not publish his proof. During Stalin's terror in the 1930's he was arrested and perished in GULag \cite{Gurin}. But a bright student Gurin who was in the class recollected the proof and published it in the {\em Americal Mathematical Monthly} as a tribute to his illustrious teacher in 1996 \cite{Gurin}. Gurin presented his intentions as follows: 
\begin{quote}
{\em To establish the power series expansions of the sine and cosine functions using only the following elementary results:
\begin{enumerate}
\item
$\displaystyle \lim_{x\rightarrow 0} \frac{sin x}{x}=1$
\item
$\displaystyle \lim_{n\rightarrow \infty}\frac{x^n}{n!}=0$
\item
Results from elementary geometry and trigonometry
\end{enumerate}
The proof must be strictly geometric and must uncover the geometric
meaning of every term in the series.}
\end{quote}

The proofs are beautiful, but requires some combinatorial arguments which may be hard even for bright students. The paper by Gurin contains the main ideas but little details.   

This paper is written with the explicit intention of bringing once again the ideas of Chaikovsky to the attention of the mathematics teachers and students and help them enhance their intuitive understanding the meaning of the mysterious power series expansions of the sine and cosine functions. The paper shows, following the footsteps of Chaikovsky, but not without using calculus, how the various terms in the power series expansions pop up geometrically.

Here is a summary of the paper. After presenting the definition of an involute in Section 2, we present the main result of the paper in Section 3. In Section 4 we present the form of equations of the involute and in Section 5 we consider the code written in the computer algebra system {\em Maxima} to compute the involutes. The results of the applications of this code are given in Section 6. Some general expressions are stated and proved in Section 7.
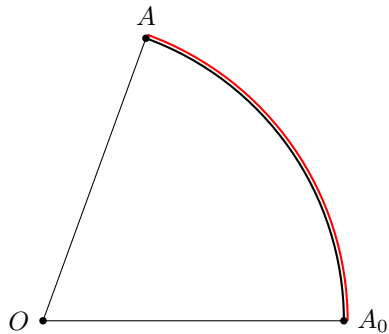
\begin{figure}[!h]
{\centering
\begin{tikzpicture}
\draw (0,0) -- (\ms,0);
\draw (0,0) -- (70: \ms);
\draw [thick, black] (\ms,0) arc (0:70:\ms);
\draw [thick, red] (\ms+0.05,0) arc (0:70:\ms + 0.05);
\draw \DrawControlLeft{(0,0)}{black}{$O$}; 
\draw \DrawControlAbove{(70:\ms)}{black}{$A$};
\draw \DrawControlRight{(\ms,0)}{black}{$A_0$};
\end{tikzpicture}
\caption{Arc of a circle (in black) and attached string (in red)}\label{Fig1}
}
\end{figure}
\begin{figure}[!h]
{\centering
\begin{tikzpicture}
\draw (0,0) -- (\ms,0);
\draw (0,0) -- (\ma r: \ms);
\draw [thick, black] (\ms,0) arc (0:\ma r:\ms);
\draw [very thick, red] (\ms,0) arc (0 r:\mc r:\ms);
\draw \DrawControlLeft{(0,0)}{black}{$O$}; 
\draw \DrawControlAbove{(\ma r:\ms)}{black}{$A$};
\draw \DrawControlRight{(\ms,0)}{black}{$A_0$};
\draw \DrawControlRight{({\ms*cos(\mc r)}, {\ms*sin(\mc r)})}{black}{$P$};
\draw[very thick, domain=0:\ma - \mc , smooth, variable = \t] 
plot 
(
{ \ms*(sin(\t r + \mb r) - \t*cos(\t r + \mb r))  }, 
{ \ms*(cos(\t r + \mb r) + \t*sin(\t r + \mb r))  } 
);
\draw  \DrawControlRight{(
{ \ms*(sin(\ma r -\mc r + \mb r) - (\ma -\mc)*cos(\ma r -\mc r + \mb r))  }, 
{ \ms*(cos(\ma r -\mc r + \mb r) + (\ma -\mc)*sin(\ma r -\mc r + \mb r))  } 
)}{black}{$P_1$};
\draw [very thick, red] (
{ (\ms)*(sin(\ma r -\mc r + \mb r) - (\ma -\mc)*cos(\ma r -\mc r + \mb r))  }, 
{ (\ms)*(cos(\ma r -\mc r + \mb r) + (\ma -\mc)*sin(\ma r -\mc r + \mb r))} ) -- ({(\ms+0.05)*cos(\mc r)}, {(\ms+0.05)*sin(\mc r)});
\end{tikzpicture}
\caption{Attached string unwinding from the arc, from the free end $A$}\label{Fig2}
}
\end{figure}
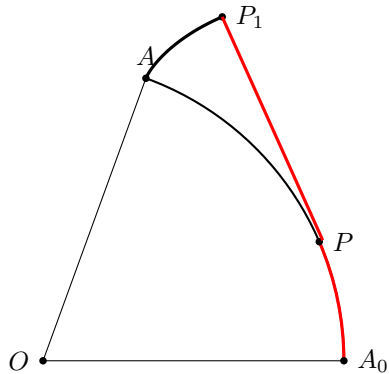
\begin{figure}[!h]
{\centering
\begin{tikzpicture}
\draw (0,0) -- (\ms,0);
\draw (0,0) -- (\ma r: \ms);
\draw [thick, black] (\ms,0) arc (0:\ma r:\ms);
\draw \DrawControlLeft{(0,0)}{black}{$O$}; 
\draw \DrawControlAbove{(\ma r:\ms)}{black}{$A$};
\draw \DrawControlRight{(\ms,0)}{black}{$A_0$};
\draw[thick, domain=0:\ma, smooth, variable = \t] 
plot 
(
{ \ms*(sin(\t r + \mb r) - \t*cos(\t r + \mb r))  }, 
{ \ms*(cos(\t r + \mb r) + \t*sin(\t r + \mb r))  } 
);
\draw  \DrawControlRight{(
{ \ms*(sin(\ma r  + \mb r) - (\ma)*cos(\ma r  + \mb r))  }, 
{ \ms*(cos(\ma r  + \mb r) + (\ma)*sin(\ma r  + \mb r))  } 
)}{black}{$A_1$};
\draw [thick, red] (
{ (\ms)*(sin(\ma r + \mb r) - (\ma)*cos(\ma r + \mb r))  }, 
{ (\ms)*(cos(\ma r + \mb r) + (\ma)*sin(\ma r + \mb r))} ) -- ({(\ms)*cos(0 r)}, {(\ms)*sin(0 r)});
\end{tikzpicture}
\caption{The curve $AA_1$ is the involute of the arc $AA_0$}\label{Fig3}
}
\end{figure}
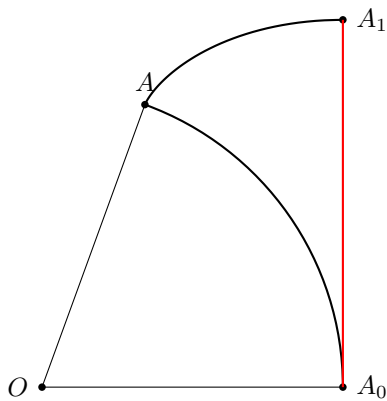
\newcommand\msc{0.5}
\begin{figure}[!b]
{\centering
\begin{tikzpicture}
\draw (0,0) -- (\ma r: \msc);
\draw [thick, black] (\msc,0) arc (0:6.2831853 r:\msc);
\draw \DrawControlLeft{(0,0)}{black}{$O$}; 
\draw \DrawControlAbove{(\ma r:\msc)}{black}{$A$};
%
\draw[thick, domain=0:6.2831853, smooth, variable = \t, blue] 
plot 
(
{ \msc*(sin(\t r + \mb r) - \t*cos(\t r + \mb r))  }, 
{ \msc*(cos(\t r + \mb r) + \t*sin(\t r + \mb r))  } 
);
\draw  \DrawControlAbove{(
{ \msc*(sin(6.2831853 r  + \mb r) - (6.2831853)*cos(6.2831853 r  + \mb r))  }, 
{ \msc*(cos(6.2831853 r  + \mb r) + (6.2831853)*sin(6.2831853 r  + \mb r))  } 
)}{black}{$A_1$};
\draw [thick, red] 
(\ma r:\msc)  -- 
( 
{ (\msc)*(sin(6.2831853 r + \mb r) - (6.2831853)*cos(6.2831853 r + \mb r))  }, 
{ (\msc)*(cos(6.2831853 r + \mb r) + (6.2831853)*sin(6.2831853 r + \mb r))});
\end{tikzpicture}
\caption{The curve $AA_1$ (shown in blue color) is an involute of the circle with center $O$}\label{Fig6}
}
\end{figure}
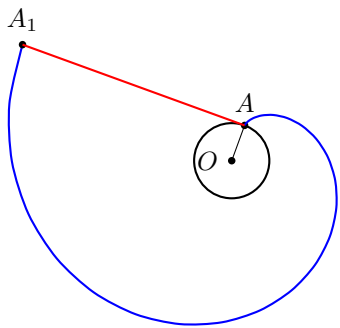
\section{Involutes}
In Chaikovsky's scheme of things, the important geometrical concept connected with the series expansions is that of an involute. 

The concept of an involute of a curve was introduced by  in 1673 by Christiaan Huygens, a Dutch mathematician and physicist. Huygenes introduced the concept in the context of his studies of the motion of a simple pendulum which itself was part of his  efforts to develop accurate time keepers. 

{\em The {\bfseries involute} of a curve $C$ is defined as a curve that is obtained by attaching an imaginary string to the curve $C$ and unwinding it tautly from the  curve $C$; the locus of the free end of this attached taut string is known as an involute of the curve $C$.} 

A given curve has an infinite number of involutes depending on where the initial position of the free end of the string lies on the curve. 

To illustrate let us see how the involute of an arc of a circle is defined. Consider an arc $AA_0$ of a circle and a string be tautly attached to the arc. Figure \ref{Fig1} shows the arc and the string (string is shown in red color).
Now the string is unwound tautly from the end $A$. Figure \ref{Fig2} shows an intermediate position of the string where the length of the line segment $PP_1$ is equal to the length of the the arc $AP$. Also the line segment $PP_1$ is tangent to the arc at $P$.
The  involute of the arc $AA_0$, which is the curve $AA_1$, is shown in Figure \ref{Fig3}. The involute of a full circle is shown in Figure \ref{Fig6}.
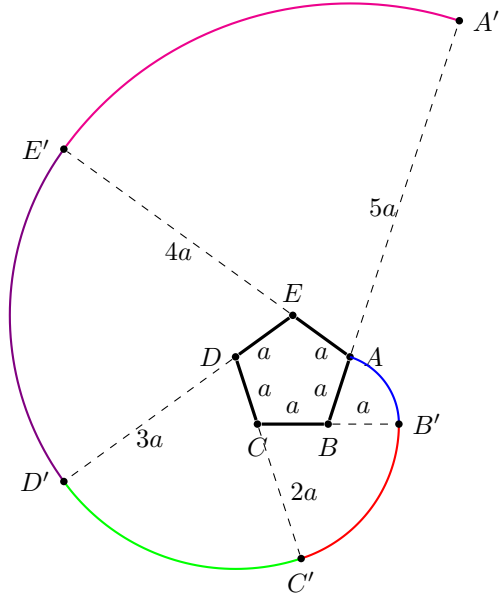
\begin{figure}
{
\centering
\begin{tikzpicture}[scale=0.4]
%
\draw node[black, circle,fill=black, inner sep=1pt,label={[black] right:{\normalsize $A$}}] (A) at (18:2) {};
\draw node[black, circle,fill=black, inner sep=1pt,label={[black]  above:{\normalsize $E$}}] (E) at (90:2) {};
\draw node[black, circle,fill=black, inner sep=1pt,label={[black] left:{\normalsize $D$}}] (D) at (162:2) {};
\draw node[black, circle,fill=black, inner sep=1pt,label={[black] below:{\normalsize $C$}}] (C) at (234:2) {};
\draw node[black, circle,fill=black, inner sep=1pt,label={[black] below:{\normalsize $B$}}] (B) at (306:2) {};
\path[draw, very thick] (A)--(B)--(C)--(D)--(E)--(A);
\draw[domain=72:0, smooth, variable = \t, color = blue, thick] 
plot 
(
{ 2*cos(306) + \side*cos(\t)  }, 
{ 2*sin(306) + \side*sin(\t)  } 
);
\draw[domain=0:-72, smooth, variable = \t, color = red, thick] 
plot 
(
{ 2*cos(234) + 2*\side*cos(\t)  }, 
{ 2*sin(234) + 2*\side*sin(\t)  } 
);
\draw[domain=-72:-144, smooth, variable = \t, color = green, thick] 
plot 
(
{ 2*cos(162) + 3*\side*cos(\t)  }, 
{ 2*sin(162) + 3*\side*sin(\t)  } 
);
\draw[domain=-144:-216, smooth, variable = \t, color = violet, thick] 
plot 
(
{ 2*cos(90) + 4*\side*cos(\t)  }, 
{ 2*sin(90) + 4*\side*sin(\t)  } 
);
\draw[domain=-216:-288, smooth, variable = \t, color = magenta, thick] 
plot 
(
{ 2*cos(18) + 5*\side*cos(\t)  }, 
{ 2*sin(18) + 5*\side*sin(\t)  } 
);
\draw node[black, circle,fill=black, inner sep=1pt,label={[black] right:{\normalsize $B^\prime$}}] (B1) at (
{ 2*cos(306) + \side*cos(0)  }, 
{ 2*sin(306) + \side*sin(0)  } 
) {};
\draw node[black, circle,fill=black, inner sep=1pt,label={[black] below:{\normalsize $C^\prime$}}] (C1) at 
(
{ 2*cos(234) + 2*\side*cos(-72)  }, 
{ 2*sin(234) + 2*\side*sin(-72)  } 
) {};
\draw node[black, circle,fill=black, inner sep=1pt,label={[black] left:{\normalsize $D^\prime$}}] (D1) at 
(
{ 2*cos(162) + 3*\side*cos(-144)  }, 
{ 2*sin(162) + 3*\side*sin(-144)  } 
) {};
\draw node[black, circle,fill=black, inner sep=1pt,label={[black] left:{\normalsize $E^\prime$}}] (E1) at 
(
{ 2*cos(90) + 4*\side*cos(-216)  }, 
{ 2*sin(90) + 4*\side*sin(-216)  } 
) {};
\draw node[black, circle,fill=black, inner sep=1pt,label={[black] right:{\normalsize $A^\prime$}}] (A1) at 
(
{ 2*cos(18) + 5*\side*cos(-288)  }, 
{ 2*sin(18) + 5*\side*sin(-288)  } 
) {};
\draw[dashed] (B)--(B1);
\draw[dashed] (C)--(C1);
\draw[dashed] (D)--(D1);
\draw[dashed] (E) -- (E1);
\draw[dashed] (A) -- (A1);
\draw node[black, inner sep = 0.25pt, label={[black] above:{\normalsize $a$}}] at (270:1.618) {};
\draw node[black, inner sep = 0.25pt, label={[black] left:{\normalsize $a$}}] at (342:1.618) {};
\draw node[black, inner sep = 0.25pt, label={[black] below:{\normalsize $a$}}] at (54:1.618) {};
\draw node[black, inner sep = 0.25pt, label={[black] below:{\normalsize $a$}}] at (126:1.618) {};
\draw node[black, inner sep = 0.25pt, label={[black] right:{\normalsize $a$}}] at (198:1.618) {};
\draw node[black, inner sep = 0.5pt, label={[black] above:{\normalsize $a$}}] at (
{ 2*cos(306) + 0.5*\side*cos(0)  }, 
{ 2*sin(306) + 0.5*\side*sin(0)  } 
) {};
\draw node[black, inner sep = 0.5pt, label={[black] right:{\normalsize $2a$}}] at 
(
{ 2*cos(234) + \side*cos(-72)  }, 
{ 2*sin(234) + \side*sin(-72)  } 
)
 {};
\draw node[black, inner sep = 0.5pt, label={[black] below:{\normalsize $3a$}}] at 
(
{ 2*cos(162) + 1.5*\side*cos(-144)  }, 
{ 2*sin(162) + 1.5*\side*sin(-144)  } 
)
 {};
\draw node[black, inner sep = 0.5pt, label={[black] below:{\normalsize $4a$}}] at (
{ 2*cos(90) + 2*\side*cos(-216)  }, 
{ 2*sin(90) + 2*\side*sin(-216)  } 
) {};
\draw node[black, inner sep = 0.5pt, label={[black] left:{\normalsize $5a$}}] at (
{ 2*cos(0) + 2.5*\side*cos(-288)  }, 
{ 2*sin(0) + 2.5*\side*sin(-288)  } 
) {};
\end{tikzpicture}
\caption{Involute of the regular pentagon $ABCDE$ is the curve $AB^\prime C^\prime D^\prime E^\prime A^\prime$}\label{Figure6}
}
\end{figure}
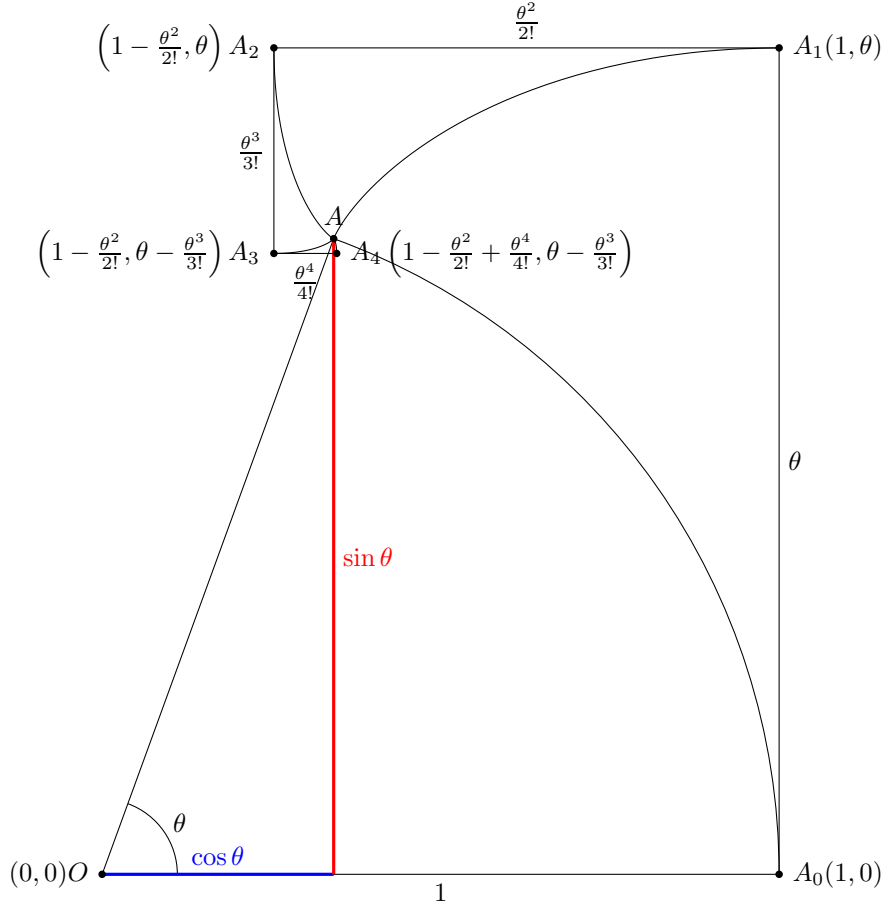
\begin{figure}[!t]
{\centering
\begin{tikzpicture}
\draw (0,0) -- (\mr,0);
\draw (0,0) -- (\ma r: \mr);
\draw [black] (\mr,0) arc (0:\ma r:\mr);
\draw [black] (1,0) arc (0:\ma r:1);
\draw [very thick, blue] (0,0) -- ( {\mr*cos(\ma r)},0);
\draw [very thick, red] ({\mr*cos(\ma r)}, {\mr*sin(\ma r)}) -- ( {\mr*cos(\ma r)},0);
\draw \DrawControlLeft{(0,0)}{black}{$(0,0)O$}; 
\draw \DrawControlAbove{(\ma r:\mr)}{black}{$A$};
\draw \DrawControlRight{(\mr,0)}{black}{$A_0(1,0)$};
\node[right]  (AA) at ({cos(0.61086 r)}, {1.2*sin(0.61086 r)}) {$\theta$};
\node[right, red] (BB) at ({\mr*cos(\ma r)}, {0.5*\mr*sin(\ma r)}) {$\sin \theta$};
\node[above, blue] (CC) at ({0.5*\mr*cos(\ma r)}, 0) {$\cos \theta$};
\draw[domain=0:\ma, smooth, variable = \t] 
plot 
(
{ \mr*(sin(\t r + \mb r) - \t*cos(\t r + \mb r))  }, 
{ \mr*(cos(\t r + \mb r) + \t*sin(\t r + \mb r))  } 
);
\draw  \DrawControlRight{(
{ \mr*(sin(\ma r  + \mb r) - (\ma)*cos(\ma r  + \mb r))  }, 
{ \mr*(cos(\ma r  + \mb r) + (\ma)*sin(\ma r  + \mb r))  } 
)}{black}{$A_1(1,\theta)$};
\draw [black] (
{ (\mr)*(sin(\ma r + \mb r) - (\ma)*cos(\ma r + \mb r))  }, 
{ (\mr)*(cos(\ma r + \mb r) + (\ma)*sin(\ma r + \mb r))} ) -- ({(\mr)*cos(0 r)}, {(\mr)*sin(0 r)});
\draw[domain=0:\ma, smooth, variable = \t] 
plot 
(
{  \mr*(sin(\t r + \mb r) - \t*cos(\t r + \mb r) - (1/2)*\t*\t*sin(\t r + \mb r)) }, 
{  \mr*(cos(\t r + \mb r) + \t*sin(\t r + \mb r) - (1/2)*\t*\t*cos(\t r + \mb r)) }
);
\draw  \DrawControlLeft{(
{  \mr*(sin(\ma r + \mb r)- \ma*cos(\ma r + \mb r)-(1/2)*\ma*\ma*sin(\ma r + \mb r)) }, 
{  \mr*(cos(\ma r + \mb r)+\ma*sin(\ma r + \mb r)-(1/2)*\ma*\ma*cos(\ma r + \mb r)) }
)
}{black}{$\left(1-\frac{\theta^2}{2!}, \theta\right) A_2$};
\draw 
(
{ (\mr)*(sin(\ma r + \mb r) - (\ma)*cos(\ma r + \mb r))  }, 
{ (\mr)*(cos(\ma r + \mb r) + (\ma)*sin(\ma r + \mb r))} 
) 
--
(
{ \mr*(sin(\ma r + \mb r)- \ma*cos(\ma r + \mb r) -(1/2)*\ma*\ma*sin(\ma r + \mb r)) }, 
{  \mr*(cos(\ma r + \mb r)+\ma*sin(\ma r + \mb r) -(1/2)*\ma*\ma*cos(\ma r + \mb r)) }
);
%
%
\draw[domain=0:\ma, smooth, variable = \t] 
plot 
(
{  \mr*(   sin(\mb r+\t r)-\t*cos(\mb r+\t r)-(1/2)*\t*\t*sin(\mb r+\t r)+(1/6)*\t*\t*\t*cos(\mb r+\t r)       ) }, 
{  \mr*(   cos(\mb r+\t r)+\t*sin(\mb r+\t r)-(1/2)*\t*\t*cos(\mb r+\t r)-(1/6)*\t*\t*\t*sin(\mb r+\t r)       ) }
);
\draw  \DrawControlLeft{
(
{  \mr*(   sin(\mb r+\ma r)-\ma*cos(\mb r+\ma r)-(1/2)*\ma*\ma*sin(\mb r+\ma r)+(1/6)*\ma*\ma*\ma*cos(\mb r+\ma r)       ) }, 
{  \mr*(   cos(\mb r+\ma r)+\ma*sin(\mb r+\ma r)-(1/2)*\ma*\ma*cos(\mb r+\ma r)-(1/6)*\ma*\ma*\ma*sin(\mb r+\ma r)       ) }
)
}{black}{$\left(1-\frac{\theta^2}{2!}, \theta -\frac{\theta^3}{3!}\right) A_3$};
\draw 
(
{  \mr*(   sin(\mb r+\ma r)-\ma*cos(\mb r+\ma r)-(1/2)*\ma*\ma*sin(\mb r+\ma r)+(1/6)*\ma*\ma*\ma*cos(\mb r+\ma r)       ) }, 
{  \mr*(   cos(\mb r+\ma r)+\ma*sin(\mb r+\ma r)-(1/2)*\ma*\ma*cos(\mb r+\ma r)-(1/6)*\ma*\ma*\ma*sin(\mb r+\ma r)       ) }
) 
--
(
{ \mr*(sin(\ma r + \mb r)- \ma*cos(\ma r + \mb r) -(1/2)*\ma*\ma*sin(\ma r + \mb r)) }, 
{  \mr*(cos(\ma r + \mb r)+\ma*sin(\ma r + \mb r) -(1/2)*\ma*\ma*cos(\ma r + \mb r)) }
);
%
%
\draw[domain=0:\ma, smooth, variable = \t] 
plot 
(
{  \mr*(   sin(\mb r+\t r)-\t*cos(\mb r+\t r)-(1/2)*\t*\t*sin(\mb r+\t r)+(1/6)*\t*\t*\t*cos(\mb r+\t r) + (1/24)*\t*\t*\t*\t*sin(\mb r + \t r)      ) }, 
{  \mr*(   cos(\mb r+\t r)+\t*sin(\mb r+\t r)-(1/2)*\t*\t*cos(\mb r+\t r)-(1/6)*\t*\t*\t*sin(\mb r+\t r) + (1/24)*\t*\t*\t*\t*cos(\mb r + \t r)     ) }
);
\draw  \DrawControlRight{
(
{  \mr*(   sin(\mb r+\ma r)-\ma*cos(\mb r+\ma r)-(1/2)*\ma*\ma*sin(\mb r+\ma r)+(1/6)*\ma*\ma*\ma*cos(\mb r+\ma r) + (1/24)*\ma*\ma*\ma*\ma*sin(\mb r + \ma r)            ) }, 
{  \mr*(   cos(\mb r+\ma r)+\ma*sin(\mb r+\ma r)-(1/2)*\ma*\ma*cos(\mb r+\ma r)-(1/6)*\ma*\ma*\ma*sin(\mb r+\ma r)  + (1/24)*\ma*\ma*\ma*\ma*cos(\mb r + \ma r)         ) }
)
}{black}{$A_4\left(1-\frac{\theta^2}{2!}+\frac{\theta^4}{4!}, \theta-\frac{\theta^3}{3!}\right)$};
\draw 
(
{  \mr*(   sin(\mb r+\ma r)-\ma*cos(\mb r+\ma r)-(1/2)*\ma*\ma*sin(\mb r+\ma r)+(1/6)*\ma*\ma*\ma*cos(\mb r+\ma r)       ) }, 
{  \mr*(   cos(\mb r+\ma r)+\ma*sin(\mb r+\ma r)-(1/2)*\ma*\ma*cos(\mb r+\ma r)-(1/6)*\ma*\ma*\ma*sin(\mb r+\ma r)       ) }
) 
--
(
{  \mr*(   sin(\mb r+\ma r)-\ma*cos(\mb r+\ma r)-(1/2)*\ma*\ma*sin(\mb r+\ma r)+(1/6)*\ma*\ma*\ma*cos(\mb r+\ma r) + (1/24)*\ma*\ma*\ma*\ma*sin(\mb r + \ma r)            ) }, 
{  \mr*(   cos(\mb r+\ma r)+\ma*sin(\mb r+\ma r)-(1/2)*\ma*\ma*cos(\mb r+\ma r)-(1/6)*\ma*\ma*\ma*sin(\mb r+\ma r)  + (1/24)*\ma*\ma*\ma*\ma*cos(\mb r + \ma r)         ) }
);
%
%
\node[below] (HH) at (0.5*\mr, 0) {$1$};
\node[right] (DD) at (\mr, \mr*0.5*\ma) {$\theta$};
\node[above] (EE) at 
( 
{  \mr*(sin(\ma r + \mb r) - \ma*cos(\ma r + \mb r) - (1/4)*\ma*\ma*sin(\ma r + \mb r)) }, 
{  \mr*(cos(\ma r + \mb r) + \ma*sin(\ma r + \mb r) - (1/2)*\ma*\ma*cos(\ma r + \mb r)) }
)
{$\frac{\theta^2}{2!}$};
\node[left] (FF) at
(
{  \mr*(sin(\ma r + \mb r)- \ma*cos(\ma r + \mb r)-(1/2)*\ma*\ma*sin(\ma r + \mb r)) },
{  \mr*(   cos(\mb r+\ma r)+\ma*sin(\mb r+\ma r)-(1/2)*\ma*\ma*cos(\mb r+\ma r)-(1/12)*\ma*\ma*\ma*sin(\mb r+\ma r)       ) }
)
{$\frac{\theta^3}{3!}$};
\node[below] (GG) at
(
{  \mr*(   sin(\mb r+\ma r)-\ma*cos(\mb r+\ma r)-(1/2)*\ma*\ma*sin(\mb r+\ma r)+(1/6)*\ma*\ma*\ma*cos(\mb r+\ma r) + (1/48)*\ma*\ma*\ma*\ma*sin(\mb r + \ma r)         ) }, 
{  \mr*(   cos(\mb r+\ma r)+\ma*sin(\mb r+\ma r)-(1/2)*\ma*\ma*cos(\mb r+\ma r)-(1/6)*\ma*\ma*\ma*sin(\mb r+\ma r) + (1/24)*\ma*\ma*\ma*\ma*cos(\mb r + \ma r)         ) }
)
{$\frac{\theta^4}{4!}$};
\end{tikzpicture}
\caption{Where do the terms of the power series expansions of the sine and cosine functions come from}\label{Fig4}
}
\end{figure}
%
%
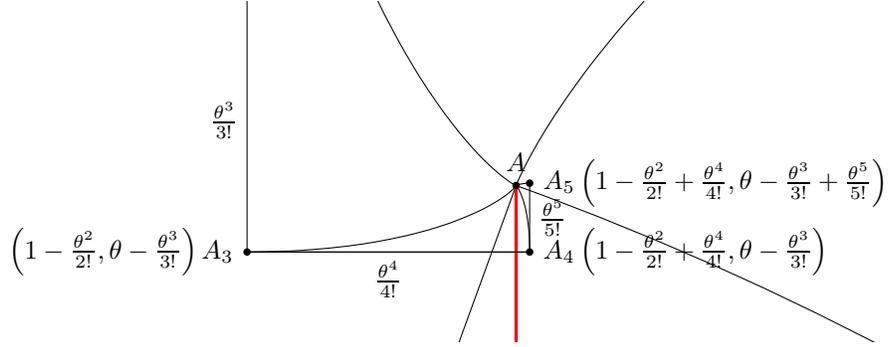
\begin{figure}
%
{
\centering
\begin{tikzpicture}[scale=4.5]
\clip (1.5,7.995) rectangle (4.3, 9); 
\draw (0,0) -- (\mr,0);
\draw (0,0) -- (\ma r: \mr);
\draw [black] (\mr,0) arc (0:\ma r:\mr);
\draw [black] (1,0) arc (0:\ma r:1);
\draw [very thick, blue] (0,0) -- ( {\mr*cos(\ma r)},0);
\draw [very thick, red] ({\mr*cos(\ma r)}, {\mr*sin(\ma r)}) -- ( {\mr*cos(\ma r)},0);
\draw \DrawControlLeft{(0,0)}{black}{$(0,0)O$}; 
\draw \DrawControlAbove{(\ma r:\mr)}{black}{$A$};
\draw \DrawControlRight{(\mr,0)}{black}{$A_0(1,0)$};
\node[right]  (AA) at ({cos(0.61086 r)}, {1.2*sin(0.61086 r)}) {$\theta$};
\node[right, red] (BB) at ({\mr*cos(\ma r)}, {0.5*\mr*sin(\ma r)}) {$\sin \theta$};
\node[above, blue] (CC) at ({0.5*\mr*cos(\ma r)}, 0) {$\cos \theta$};
\draw[domain=0:\ma, smooth, variable = \t] 
plot 
(
{ \mr*(sin(\t r + \mb r) - \t*cos(\t r + \mb r))  }, 
{ \mr*(cos(\t r + \mb r) + \t*sin(\t r + \mb r))  } 
);
\draw  \DrawControlRight{(
{ \mr*(sin(\ma r  + \mb r) - (\ma)*cos(\ma r  + \mb r))  }, 
{ \mr*(cos(\ma r  + \mb r) + (\ma)*sin(\ma r  + \mb r))  } 
)}{black}{$A_1(1,\theta)$};
\draw [black] (
{ (\mr)*(sin(\ma r + \mb r) - (\ma)*cos(\ma r + \mb r))  }, 
{ (\mr)*(cos(\ma r + \mb r) + (\ma)*sin(\ma r + \mb r))} ) -- ({(\mr)*cos(0 r)}, {(\mr)*sin(0 r)});
\draw[domain=0:\ma, smooth, variable = \t] 
plot 
(
{  \mr*(sin(\t r + \mb r) - \t*cos(\t r + \mb r) - (1/2)*\t*\t*sin(\t r + \mb r)) }, 
{  \mr*(cos(\t r + \mb r) + \t*sin(\t r + \mb r) - (1/2)*\t*\t*cos(\t r + \mb r)) }
);
\draw  \DrawControlLeft{(
{  \mr*(sin(\ma r + \mb r)- \ma*cos(\ma r + \mb r)-(1/2)*\ma*\ma*sin(\ma r + \mb r)) }, 
{  \mr*(cos(\ma r + \mb r)+\ma*sin(\ma r + \mb r)-(1/2)*\ma*\ma*cos(\ma r + \mb r)) }
)
}{black}{$\left(1-\frac{\theta^2}{2!}, \theta\right) A_2$};
\draw 
(
{ (\mr)*(sin(\ma r + \mb r) - (\ma)*cos(\ma r + \mb r))  }, 
{ (\mr)*(cos(\ma r + \mb r) + (\ma)*sin(\ma r + \mb r))} 
) 
--
(
{ \mr*(sin(\ma r + \mb r)- \ma*cos(\ma r + \mb r) -(1/2)*\ma*\ma*sin(\ma r + \mb r)) }, 
{  \mr*(cos(\ma r + \mb r)+\ma*sin(\ma r + \mb r) -(1/2)*\ma*\ma*cos(\ma r + \mb r)) }
);
%
%
\draw[domain=0:\ma, smooth, variable = \t] 
plot 
(
{  \mr*(   sin(\mb r+\t r)-\t*cos(\mb r+\t r)-(1/2)*\t*\t*sin(\mb r+\t r)+(1/6)*\t*\t*\t*cos(\mb r+\t r)       ) }, 
{  \mr*(   cos(\mb r+\t r)+\t*sin(\mb r+\t r)-(1/2)*\t*\t*cos(\mb r+\t r)-(1/6)*\t*\t*\t*sin(\mb r+\t r)       ) }
);
\draw  \DrawControlLeft{
(
{  \mr*(   sin(\mb r+\ma r)-\ma*cos(\mb r+\ma r)-(1/2)*\ma*\ma*sin(\mb r+\ma r)+(1/6)*\ma*\ma*\ma*cos(\mb r+\ma r)       ) }, 
{  \mr*(   cos(\mb r+\ma r)+\ma*sin(\mb r+\ma r)-(1/2)*\ma*\ma*cos(\mb r+\ma r)-(1/6)*\ma*\ma*\ma*sin(\mb r+\ma r)       ) }
)
}{black}{$\left(1-\frac{\theta^2}{2!}, \theta -\frac{\theta^3}{3!}\right) A_3$};
\draw 
(
{  \mr*(   sin(\mb r+\ma r)-\ma*cos(\mb r+\ma r)-(1/2)*\ma*\ma*sin(\mb r+\ma r)+(1/6)*\ma*\ma*\ma*cos(\mb r+\ma r)       ) }, 
{  \mr*(   cos(\mb r+\ma r)+\ma*sin(\mb r+\ma r)-(1/2)*\ma*\ma*cos(\mb r+\ma r)-(1/6)*\ma*\ma*\ma*sin(\mb r+\ma r)       ) }
) 
--
(
{ \mr*(sin(\ma r + \mb r)- \ma*cos(\ma r + \mb r) -(1/2)*\ma*\ma*sin(\ma r + \mb r)) }, 
{  \mr*(cos(\ma r + \mb r)+\ma*sin(\ma r + \mb r) -(1/2)*\ma*\ma*cos(\ma r + \mb r)) }
);
%
%
\draw[domain=0:\ma, smooth, variable = \t] 
plot 
(
{  \mr*(   sin(\mb r+\t r)-\t*cos(\mb r+\t r)-(1/2)*\t*\t*sin(\mb r+\t r)+(1/6)*\t*\t*\t*cos(\mb r+\t r) + (1/24)*\t*\t*\t*\t*sin(\mb r + \t r)      ) }, 
{  \mr*(   cos(\mb r+\t r)+\t*sin(\mb r+\t r)-(1/2)*\t*\t*cos(\mb r+\t r)-(1/6)*\t*\t*\t*sin(\mb r+\t r) + (1/24)*\t*\t*\t*\t*cos(\mb r + \t r)     ) }
);
\draw  \DrawControlRight{
(
{  \mr*(   sin(\mb r+\ma r)-\ma*cos(\mb r+\ma r)-(1/2)*\ma*\ma*sin(\mb r+\ma r)+(1/6)*\ma*\ma*\ma*cos(\mb r+\ma r) + (1/24)*\ma*\ma*\ma*\ma*sin(\mb r + \ma r)            ) }, 
{  \mr*(   cos(\mb r+\ma r)+\ma*sin(\mb r+\ma r)-(1/2)*\ma*\ma*cos(\mb r+\ma r)-(1/6)*\ma*\ma*\ma*sin(\mb r+\ma r)  + (1/24)*\ma*\ma*\ma*\ma*cos(\mb r + \ma r)         ) }
)
}{black}{$A_4\left(1-\frac{\theta^2}{2!}+\frac{\theta^4}{4!}, \theta-\frac{\theta^3}{3!}\right)$};
\draw 
(
{  \mr*(   sin(\mb r+\ma r)-\ma*cos(\mb r+\ma r)-(1/2)*\ma*\ma*sin(\mb r+\ma r)+(1/6)*\ma*\ma*\ma*cos(\mb r+\ma r)       ) }, 
{  \mr*(   cos(\mb r+\ma r)+\ma*sin(\mb r+\ma r)-(1/2)*\ma*\ma*cos(\mb r+\ma r)-(1/6)*\ma*\ma*\ma*sin(\mb r+\ma r)       ) }
) 
--
(
{  \mr*(   sin(\mb r+\ma r)-\ma*cos(\mb r+\ma r)-(1/2)*\ma*\ma*sin(\mb r+\ma r)+(1/6)*\ma*\ma*\ma*cos(\mb r+\ma r) + (1/24)*\ma*\ma*\ma*\ma*sin(\mb r + \ma r)            ) }, 
{  \mr*(   cos(\mb r+\ma r)+\ma*sin(\mb r+\ma r)-(1/2)*\ma*\ma*cos(\mb r+\ma r)-(1/6)*\ma*\ma*\ma*sin(\mb r+\ma r)  + (1/24)*\ma*\ma*\ma*\ma*cos(\mb r + \ma r)         ) }
);
%
%
\draw[domain=0:\ma, smooth, variable = \t] 
plot 
(
{  \mr*(   sin(\mb r+\t r)-\t*cos(\mb r+\t r)-(1/2)*\t*\t*sin(\mb r+\t  r)+(1/6)*\t*\t*\t*cos(\mb r+\t r) + (1/24)*\t*\t*\t*\t*sin(\mb r + \t r)       - (1/120)*\t*\t*\t*\t*\t*cos(\mb r+\t r))}, 
{  \mr*(   cos(\mb r+\t r)+\t*sin(\mb r+\t r)-(1/2)*\t*\t*cos(\mb r+\t r)-(1/6)*\t*\t*\t*sin(\mb r+\t r) + (1/24)*\t*\t*\t*\t*cos(\mb r + \t r)      + (1/120)*\t*\t*\t*\t*\t*sin(\mb r+\t r))}
);
\draw  \DrawControlRight{
(
{  \mr*(   sin(\mb r+\ma r)-\ma*cos(\mb r+\ma r)-(1/2)*\ma*\ma*sin(\mb r+\ma r)+(1/6)*\ma*\ma*\ma*cos(\mb r+\ma r) + (1/24)*\ma*\ma*\ma*\ma*sin(\mb r + \ma r) - (1/120)*\ma*\ma*\ma*\ma*\ma*cos(\mb r+\ma r)           ) }, 
{  \mr*(   cos(\mb r+\ma r)+\ma*sin(\mb r+\ma r)-(1/2)*\ma*\ma*cos(\mb r+\ma r)-(1/6)*\ma*\ma*\ma*sin(\mb r+\ma r)  + (1/24)*\ma*\ma*\ma*\ma*cos(\mb r + \ma r) + (1/120)*\ma*\ma*\ma*\ma*\ma*sin(\mb r+\ma r)        ) }
)
}{black}{$A_5\left(1-\frac{\theta^2}{2!}+\frac{\theta^4}{4!}, \theta-\frac{\theta^3}{3!} + \frac{\theta^5}{5!}\right)$};
\draw 
(
{  \mr*(   sin(\mb r+\ma r)-\ma*cos(\mb r+\ma r)-(1/2)*\ma*\ma*sin(\mb r+\ma r)+(1/6)*\ma*\ma*\ma*cos(\mb r+\ma r)       ) }, 
{  \mr*(   cos(\mb r+\ma r)+\ma*sin(\mb r+\ma r)-(1/2)*\ma*\ma*cos(\mb r+\ma r)-(1/6)*\ma*\ma*\ma*sin(\mb r+\ma r)       ) }
) 
--
(
{  \mr*(   sin(\mb r+\ma r)-\ma*cos(\mb r+\ma r)-(1/2)*\ma*\ma*sin(\mb r+\ma r)+(1/6)*\ma*\ma*\ma*cos(\mb r+\ma r) + (1/24)*\ma*\ma*\ma*\ma*sin(\mb r + \ma r)            ) }, 
{  \mr*(   cos(\mb r+\ma r)+\ma*sin(\mb r+\ma r)-(1/2)*\ma*\ma*cos(\mb r+\ma r)-(1/6)*\ma*\ma*\ma*sin(\mb r+\ma r)  + (1/24)*\ma*\ma*\ma*\ma*cos(\mb r + \ma r)         ) }
);
\draw 
(
{  \mr*(   sin(\mb r+\ma r)-\ma*cos(\mb r+\ma r)-(1/2)*\ma*\ma*sin(\mb r+\ma r)+(1/6)*\ma*\ma*\ma*cos(\mb r+\ma r) + (1/24)*\ma*\ma*\ma*\ma*sin(\mb r + \ma r) - (1/120)*\ma*\ma*\ma*\ma*\ma*cos(\mb r+\ma r)           ) }, 
{  \mr*(   cos(\mb r+\ma r)+\ma*sin(\mb r+\ma r)-(1/2)*\ma*\ma*cos(\mb r+\ma r)-(1/6)*\ma*\ma*\ma*sin(\mb r+\ma r)  + (1/24)*\ma*\ma*\ma*\ma*cos(\mb r + \ma r) + (1/120)*\ma*\ma*\ma*\ma*\ma*sin(\mb r+\ma r)        ) }
) 
--
(
{  \mr*(   sin(\mb r+\ma r)-\ma*cos(\mb r+\ma r)-(1/2)*\ma*\ma*sin(\mb r+\ma r)+(1/6)*\ma*\ma*\ma*cos(\mb r+\ma r) + (1/24)*\ma*\ma*\ma*\ma*sin(\mb r + \ma r)            ) }, 
{  \mr*(   cos(\mb r+\ma r)+\ma*sin(\mb r+\ma r)-(1/2)*\ma*\ma*cos(\mb r+\ma r)-(1/6)*\ma*\ma*\ma*sin(\mb r+\ma r)  + (1/24)*\ma*\ma*\ma*\ma*cos(\mb r + \ma r)         ) }
);
%
%
\node[below] (HH) at (0.5*\mr, 0) {$1$};
\node[right] (DD) at (\mr, \mr*0.5*\ma) {$\theta$};
\node[above] (EE) at 
( 
{  \mr*(sin(\ma r + \mb r) - \ma*cos(\ma r + \mb r) - (1/4)*\ma*\ma*sin(\ma r + \mb r)) }, 
{  \mr*(cos(\ma r + \mb r) + \ma*sin(\ma r + \mb r) - (1/2)*\ma*\ma*cos(\ma r + \mb r)) }
)
{$\frac{\theta^2}{2!}$};
\node[left] (FF) at
(
{  \mr*(sin(\ma r + \mb r)- \ma*cos(\ma r + \mb r)-(1/2)*\ma*\ma*sin(\ma r + \mb r)) },
{  \mr*(   cos(\mb r+\ma r)+\ma*sin(\mb r+\ma r)-(1/2)*\ma*\ma*cos(\mb r+\ma r)-(1/7)*\ma*\ma*\ma*sin(\mb r+\ma r)       ) }
)
{$\frac{\theta^3}{3!}$};
\node[below] (GG) at
(
{  \mr*(   sin(\mb r+\ma r)-\ma*cos(\mb r+\ma r)-(1/2)*\ma*\ma*sin(\mb r+\ma r)+(1/6)*\ma*\ma*\ma*cos(\mb r+\ma r) + (1/48)*\ma*\ma*\ma*\ma*sin(\mb r + \ma r)         ) }, 
{  \mr*(   cos(\mb r+\ma r)+\ma*sin(\mb r+\ma r)-(1/2)*\ma*\ma*cos(\mb r+\ma r)-(1/6)*\ma*\ma*\ma*sin(\mb r+\ma r) + (1/24)*\ma*\ma*\ma*\ma*cos(\mb r + \ma r)         ) }
)
{$\frac{\theta^4}{4!}$};
\node[right] (II) at
(
{  \mr*(   sin(\mb r+\ma r)-\ma*cos(\mb r+\ma r)-(1/2)*\ma*\ma*sin(\mb r+\ma r)+(1/6)*\ma*\ma*\ma*cos(\mb r+\ma r) + (1/24)*\ma*\ma*\ma*\ma*sin(\mb r + \ma r) - (1/240)*\ma*\ma*\ma*\ma*\ma*cos(\mb r+\ma r)        ) }, 
{  \mr*(   cos(\mb r+\ma r)+\ma*sin(\mb r+\ma r)-(1/2)*\ma*\ma*cos(\mb r+\ma r)-(1/6)*\ma*\ma*\ma*sin(\mb r+\ma r) + (1/24)*\ma*\ma*\ma*\ma*cos(\mb r + \ma r) + (1/240)*\ma*\ma*\ma*\ma*\ma*sin(\mb r+\ma r)         ) }
)
{$\frac{\theta^5}{5!}$};
\end{tikzpicture}
\caption{A zoomed in picture}\label{Fig5}
}
\end{figure}

Though not much related to the main ideas of this paper, it is interesting to note that it is indeed possible to define involutes of non-smooth curves.  For example one can construct the involutes of polygons. Figure \ref{Figure6} shows the involute of a regular pentagon. 
%
%
\section{Main result}
The data given in Figure \ref{Fig4} summarizes the main result of this paper.

Consider an arc $AA_0$ of a circle of unit radius subtending an angle of $\theta$ radians at the center $O$ of the circle. We choose $O$ as origin, the line $OA_0$ as the positive direction of the $x$-axis and an appropriate $y$-axis. With these axes, the points are $A(\cos \theta, \sin \theta)$ and $A_0(1,0)$. 

In Figure \ref{Fig4}, the curve $AA_1$ is the involute of the circular arc $AA_0$. The coordinates of $A_1$ are $(1,\theta)$ so that the length of the line segment $A_0A_1$ is $\theta$. The curve $AA_2$ is the involute of the curve $AA_1$. We show below that $A_2$ is $\left(1-\frac{\theta^2}{2!},\theta\right)$ so that 
the length of the line segment $A_1A_2$ is $\frac{\theta^2}{2!}$. The involute of the curve $AA_2$ is the curve $AA_3$ and we show that $A_3$ is $\left( 1-\frac{\theta^2}{2!}, \theta-\frac{\theta^3}{3!}\right)$, so that the length of the segment $A_3A_4$ is $\frac{\theta^3}{3!}$. We continue this process. 

The general expressions for the coordinates of the points $A_k$  can  be obtained as follows:
\begin{align*}
A_{2n-1} & \left(1-\frac{\theta^2}{2!}+\cdots+(-1)^{n-1}\frac{\theta^{2(n-1)}}{(2(n-1))!}, \theta-\frac{\theta^3}{3!}+\cdots+ (-1)^{n-1}\frac{\theta^{2n-1}}{(2n-1)!}\right)\\
A_{2n} &\left(1-\frac{\theta^2}{2!}+\cdots+(-1)^{n-1}\frac{\theta^{2n}}{(2n)!}, \theta-\frac{\theta^3}{3!}+\cdots+ (-1)^{i-1}\frac{\theta^{2n-1}}{(2n-1)!}\right)\\
\end{align*}
These expressions will be derived in the next two sections of this paper.

The coordinates of the sequence of points $A_r$ are the partial sums of the power series expansions of $\cos\theta$ and $\sin\theta$ and hence as $r\rightarrow \infty$, $A_r\rightarrow A$.

\subsection*{Main result}
\begin{itemize}
\item
The terms of the power series expansion of $\sin \theta$ are the lengths of the line segments $A_0A_1$, $A_2A_3$, $A_4A_5$, . . .
\item
The terms of the power series expansion of $\cos\theta$ are the lengths of the line segments $OA_0$, $A_1A_2$, $A_3A_4$, . . .
\end{itemize}
\section{Computation of the involutes}
Let the parametric equations of a curve be
\begin{equation}\label{curve}
x=f(t), \qquad y= g(t)\qquad a\le t\le b.
\end{equation}
The parametric equations of the involute are given by (see, for example, \cite{Oprea03})
\begin{align}
x &=f(t) - s(t) \frac{f^\prime(t)}{\sqrt{(f^\prime(t))^2 + (g^\prime(t))^2}}\label{involute1}\\
y &=g(t) - s(t) \frac{g^\prime(t)}{\sqrt{(f^\prime(t))^2 + (g^\prime(t))^2}}\label{involute2}
\end{align}
where $s(t)$ is the length of the arc of the curve given by Eq.\eqref{curve} from the initial point $t=a$ to an arbitrary point specified by the parameter $t$ and is computed using 
\begin{equation}\label{st}
s(t)=\int_a^t \sqrt{(f^\prime(u))^2 + (g^\prime(u))^2}\, du.
\end{equation}

The computations can be carried out by hand. But for the problem at hand where we have to compute the involutes of the arcs $AA_1$, $AA_2$, $AA_3$, etc. shown in Figure \ref{Fig4}, the computations are indeed daunting. Fortunately the powers of any computer algebra system can be harnessed to simplify the computations. 
\section{Code for computing involute}
The tedious computations involved in the determination of the
parametric equations of the involute of a  curve can be easily
carried out using a computer  algebra system. We give below a
procedure for this written in {\em Maxima} which is a free 
software released under the terms of the GNU General Public License \cite{Maxima}.

For a procedure written in the computer algebra system {\em Maple} the reader is referred to the relevant {\em Maple} help page written by John Oprea (see \cite{Oprea02}). (The reader is also referred to an alternative, slightly simpler, procedure given in \cite{Oprea01}). 
\subsection*{The  code}
\begin{verbatim}

/* Maxima procedure to print the parametric equations of 
   the involute of a given curve. */

   restart$ 
   assume(t>0)$ 

/* We begin by loading the parametric equations of the curve. 
   We assume that the initial value of t is 0. */
   x(t) := f(t)$ 
   y(t) := g(t)$

/* The next function simplifies some of the expressions 
   in the commands given later in the procedure. */
   norm(a,b) := sqrt(a^2 + b^2)$

/* We denote by s(t) is the length of the arc of the given 
   curve from t=0 to an arbitrary point 't' on the curve. 
   The next command defines the function s(t).*/
   s(t) := integrate(norm(diff(x(u),u),diff(y(u),u)), u,0,t)$

/* The vector (tx(t), ty(t)) is the unit tangent vector 
   to the curve at the point 't' on the curve. */
   tx(t):= diff(x(t),t)/norm(diff(x(t),t), diff(y(t),t))$
   ty(t):= diff(y(t),t)/norm(diff(x(t),t), diff(y(t),t))$

/* x = ix(t), y = iy(t) are the parametric equations of the 
   involute. The next two commands define the functions 
   ix(t) and iy(t). We may use commands like trigreduce 
   to simplify the resulting expressions.  */
   ix(t) := x(t) - s(t)*tx(t)$
   iy(t) := y(t) - s(t)*ty(t)$

/* We now print the expressions which define the parametric 
   equations of the involute. */
   ix(t);
   iy(t);

/* When we run this procedure we get the following output.*/

\end{verbatim}

\begin{math}\displaystyle
\mathrm{f}\left( t\right) -\frac{\left( \frac{d}{d\,t}\cdot \mathrm{f}\left( t\right) \right) \cdot \int_{0}^{t}\sqrt{{{\left( \frac{d}{d\,u}\cdot \mathrm{g}\left( u\right) \right) }^{2}}+{{\left( \frac{d}{d\,u}\cdot \mathrm{f}\left( u\right) \right) }^{2}}}du}{\sqrt{{{\left( \frac{d}{d\,t}\cdot \mathrm{g}\left( t\right) \right) }^{2}}+{{\left( \frac{d}{d\,t}\cdot \mathrm{f}\left( t\right) \right) }^{2}}}}
\end{math}
\begin{math}\displaystyle
\mathrm{g}\left( t\right) -\frac{\left( \frac{d}{d\,t}\cdot \mathrm{g}\left( t\right) \right) \cdot \int_{0}^{t}\sqrt{{{\left( \frac{d}{d\,u}\cdot \mathrm{g}\left( u\right) \right) }^{2}}+{{\left( \frac{d}{d\,u}\cdot \mathrm{f}\left( u\right) \right) }^{2}}}du}{\sqrt{{{\left( \frac{d}{d\,t}\cdot \mathrm{g}\left( t\right) \right) }^{2}}+{{\left( \frac{d}{d\,t}\cdot \mathrm{f}\left( t\right) \right) }^{2}}}}\mbox{}
\end{math}
\begin{verbatim}

/* Example of application of the procedure: 
   Let us print the expressions defining the parametric
   equations of the involute of the circle defined by
   x = cos t, y = sin t. */

   x(t):= cos(t); 
   y(t):= sin(t);
   ix(t);
   iy(t);

/* This produces the following output. */
\end{verbatim}

\begin{math}\displaystyle
\mathrm{x}\left( t\right) :=\mathrm{cos}\left( t\right) 
\end{math}

\begin{math}\displaystyle
\mathrm{y}\left( t\right) :=\mathrm{sin}\left( t\right) 
\end{math}

\begin{math}\displaystyle
t\cdot \mathrm{sin}\left( t\right) +\mathrm{cos}\left( t\right) \end{math}

\begin{math}\displaystyle
\mathrm{sin}\left( t\right) -t\cdot \mathrm{cos}\left( t\right) \mbox{}
\end{math}
\section{Equations of the involutes}\label{S6}
In this section, we first consider the the parametric equations of the arc $AA_0$. 
We then give the parametric equations of the  curves   $AA_1$, $AA_2$, $AA_3$ and $AA_4$ obtained by applying the {\em Maxima} procedure described above. 

Based on these computations, we present in the next section
the equations of the general involute $AA_n$.  The validity of these general expressions will be verified by induction on $n$.

To simplify the computations we set 
$$
\phi = \frac{\pi}{2}-\theta
$$
\subsubsection*{Equations of the arc $AA_0$}
The parametric equations of the arc $AA_0$ can be written as
\begin{equation}\label{AA0}
x  = \sin(t + \phi), \qquad
y  = \cos(t + \phi)\qquad 0\le t\le \theta.
\end{equation}
When $t=\theta$, we get the point $A_0(1,0)$.
\subsubsection*{Equations of the curve $AA_1$}
\begin{align}
x & = \sin(\phi+t)-t\cos(\phi+t)\label{AA1x}\\
y & = \cos(\phi+t)+t\sin(\phi+t)\label{AA1y}
\end{align}
When $t=\theta$ we get the point $A_1(1,\theta)$.
\subsubsection*{Equations of the curve $AA_2$}
\begin{align*}
x & = \sin(\phi+t)-t\cos(\phi+t)-(1/2)t^2\sin(\phi+t)\\
y & = \cos(\phi+t)+t\sin(\phi+t)-(1/2)t^2\cos(\phi+t)
\end{align*}
When $t=\theta$ we get the point $A_2\left(1-\frac{\theta^2}{2}, \theta\right)$
\subsubsection*{Equations of the curve $AA_3$}
\begin{align*}
x & = \sin(\phi+t)-t\cos(\phi+t)-(1/2)t^2\sin(\phi+t)+(1/6)t^3\cos(\phi+t)\\
y & = \cos(\phi+t)+t\sin(\phi+t)-(1/2)t^2\cos(\phi+t)-(1/6)t^3
\sin(\phi+t)
\end{align*}
When $t=\theta$ we get the point $A_3\left( 1-\frac{\theta^2}{2}, \theta -\frac{\theta^3}{6}  \right)$
\subsubsection*{Equations of the curve $AA_4$}
\begin{align*}
x & = \sin(\phi+t)-
t\cos(\phi+t)-(1/2)t^2\sin(\phi+t) +(1/6)t^3\cos(\phi+t)\\
& \phantom{=} +(1/24)t^4\sin(\phi+t)\\ 
y & =  \cos(\phi+t) +t\sin(\phi+t)-(1/2)t^2\cos(\phi+t) - (1/6)t^3\sin(\phi+t)\\
&\phantom{=} +(1/24)t^4\cos(\phi+t)
\end{align*}
When $t=\theta$ we get the point $A_4\left( 1-\frac{\theta^2}{2} +\frac{\theta^4}{4!}, \theta -\frac{\theta^3}{6}  \right)$
\section{General expressions}
We use the following notations:
\begin{align*}
C_n(t) & = \sum_{i=0}^{n} \frac{(-1)^{i}t^{2i}}{(2i)!},\\
S_n(t) & = \sum_{i=0}^{n-1} \frac{(-1)^{i} t^{2i+1}}{(2i+1)!}.
\end{align*}
These are the partial sums of the power series expansions of the the cosine and sine functions: $C_n$ is the sum of the first $n+1$ terms in the expansion of $\cos t$ and $S_n$ is the sum of the first $n$ terms in the expansion of $\sin t$. 
We note that 
\begin{align}
C^\prime_n (t) & = -S_n(t)\label{D1} \\
S^\prime_n (t) & = C_{n-1}(t)\label{D2}
\end{align}
\subsubsection*{Proposition}
Let the parametric equations of the involute $AA_k$ be 
\begin{align*}
x & = x_k(t), \\
y & = y_k(t).
\end{align*}
\begin{enumerate}
\item
Let $k$ be odd and $k=2n-1$ for some integer $n>0$. Then
\begin{align}
x_{k}(t)  & =
C_{n-1}(t) \sin(\phi + t) - S_{n}(t) \cos(\phi+t)\label{Ox} \\
y_{k}(t)  & = 
C_{n-1}(t) \cos(\phi + t) + S_{n}(t) \sin(\phi+t)\label{Oy}
\end{align}
\item
Let $k$ be even and $k=2n$ for some integer $n>0$. Then 
\begin{align}
x_{k}(t) & =  C_n(t) \sin(\phi + t) - S_n(t) \cos(\phi+t)\label{Ex} \\
y_{k}(t) & =  C_n(t) \cos(\phi + t) + S_n(t) \sin(\phi+t) \label{Ey}
\end{align}
\end{enumerate}
\subsection*{Proof}
We prove these expressions by induction on $k$. 
\subsubsection*{Step 1}
The result is true when $k=1$. We have already noted this in Section \ref{S6} (see Eq.\eqref{AA1x} and Eq.\eqref{AA1y}).
\subsubsection*{Step 2}
Let the result be true for some odd integer $k=2m-1\ge 1$, that is, let
\begin{align*}
x_{k}(t)  & =
C_{m-1}(t) \sin(\phi + t) - S_{m}(t) \cos(\phi+t) \\
y_{k}(t)  & = 
C_{m-1}(t) \cos(\phi + t) + S_{m}(t) \sin(\phi+t)
\end{align*}
Differentiating these expressions and using Eq.\eqref{D1} and Eq.\eqref{D2}, we get
\begin{align*}
x^\prime_k(t) & = \frac{(-1)^{m-1}t^{2m-1}}{(2m-1)!}\sin(\phi+t)\\
y^\prime_k(t) & = \frac{(-1)^{m-1}t^{2m-1}}{(2m-1)!}\cos(\phi+t)
\end{align*}
Computing the arc distance $s(t)$ using Eq.\eqref{st} we get
$$
s(t) = \frac{t^{2m}}{(2m)!}.
$$
Now using Eq.\eqref{involute1} and Eq.\eqref{involute2} and simplifying the expressions we get
\begin{align*}
x_{k+1}(t) & = C_m(t) \sin(\phi + t) - S_m(t) \cos(\phi+t) \\ 
y_{k+1}(t) & = C_m(t) \cos(\phi + t) + S_m(t) \sin(\phi+t)
\end{align*}
The expressions on the right hand side agrees with the expressions in the right hand sides of Eq.\eqref{Ex} and Eq.\eqref{Ey} with $n=m$.

Thus, if the proposition is true for some odd integer $k\ge 1$, then it is also true for next even integer $k+1$. 
\subsubsection*{Step 3}
Doing the computations as in Step 3, one can show that if the proposition is true for some even integer $k>0$ then it is also true for the next odd integer $k+1$.
\subsubsection*{Step 4}
Putting all the pieces together we conclude that the proposition si true for all integers $k\ge 1$. $\blacksquare$
\subsubsection*{Remark}
Substituting $\theta$ for $t$ in the expressions given in the proposition above we get the coordinates of the point $A_n$ as given in Section 3.
%
%

%
%
%
\end{document}